\documentclass[12pt]{article}                                               

\input{amssym.def}                                                           
\input{amssym.tex}

\begin{document}                                                             
\title{On a problem of A.~Weil}

\author{Igor ~V. ~Nikolaev}


\date{}
 \maketitle


\newtheorem{thm}{Theorem}
\newtheorem{lem}{Lemma}
\newtheorem{dfn}{Definition}
\newtheorem{rmk}{Remark}
\newtheorem{cor}{Corollary}
\newtheorem{con}{Conjecture}
\newtheorem{prp}{Proposition}


\newcommand{\N}{{\Bbb N}}
\newcommand{\F}{{\cal F}}
\newcommand{\R}{{\Bbb R}}
\newcommand{\Z}{{\Bbb Z}}
\newcommand{\C}{{\Bbb C}}

\begin{abstract}
A  topological invariant of the
geodesic laminations on a modular surface is constructed. 
The invariant has a continuous part (the tail  of
a continued fraction)  and a combinatorial part  (the singularity
data). It is shown,   that the invariant is complete, i.e. the geodesic 
lamination can be  recovered from the  invariant. The continuous part 
of the invariant has geometric meaning of a slope of lamination on 
the surface.

\vspace{7mm}


\vspace{5mm}
\end{abstract}


\section{Introduction}
Let $S$ be a connected complete hyperbolic surface. By a {\it geodesic} {\bf g}
on $S$ we understand the maximal arc consisting of the locally shortest sub-arcs;
the  geodesic {\bf g} is called {\it simple} if it has no self-crossing
points.  It is well known,  that  simple geodesic
is either (i) a closed geodesic, (ii) a non-closed spiral geodesic tending to
a closed geodesic or (iii) a non-closed geodesic,  whose limit set is a
perfect (Cantor) subset of $S$ [Casson \&  Bleiler   1988]   \cite{CaB}. 
 A  {\it geodesic lamination} is a closed subset of the surface $S$,  which is the union of 
disjoint simple geodesics {\bf g}.

The Chabauty topology (also known as the Gromov-Hausdorff topology) 
 turns the set of minimal laminations  into an important topological space
$\Lambda=\Lambda(S)$. For example, certain compactifications of the  Teichm\"uller space 
of $S$ are homeomorphic to $\Lambda$, modulo 
the laminations with more than one independent ergodic measure [Thurston 1997]   \cite{T}.  
The space $\Lambda$ is a compact Hausdorff topological space;  it has a metric measuring the
angles between the asymptotic direction (a slope) of the laminations on the surface
$S$. In the simplest case $S=T^2$ (a flat torus),  all possible slopes are exhausted by
the irrationals $\theta\in [0,2\pi[$; the latter are known as the Poincar\'e rotation numbers.

In 1936 A.~Weil asked about a generalization of the rotation numbers
to the case of the higher genus surfaces [Weil 1936]   \cite{Wei}. 
Namely, the {\it Weil problem} consists in an explicit construction
of the rotation numbers for laminations on surfaces of genus $g\ge 1$;
the numbers must satisfy all formal properties of the Poincar\'e rotation numbers.   
It was conjectured,  that the hyperbolic plane 
might be critical to a solution of the problem, {\it ibid}.
An excellent survey of  [Anosov 1995]  \cite{Ano} gives an
account of the Weil's problem after 1936.  Many important
contributions to a solution of the problem are due to  [Anosov  \& Zhuzhoma 2005]   \cite{AZ},  
 [Moeckel  1982]   \cite{Moe},  [Schwartzman 1957]   \cite{Sch},   [Series 1985]   \cite{Ser}
and others, see [Anosov 1995]  \cite{Ano}. Let us mention an influential and  spiritually close 
work of   [Artin  1924]   \cite{Art}.    
Needless to say, a solution to the Weil problem is significant
for hyperbolic geometry and low-dimensional topology.

In this note we define  a slope $\theta$  of lamination $\lambda$ on a surface $S$. 
The slope  measures  an  asymptotic direction of $\lambda$  on $S$. 
To formalize our result, we assume that $S\cong {\Bbb H}^*/G$,
where $G$ is a finite index subgroup of the modular group $SL(2,{\Bbb Z})$
and  ${\Bbb H}^*={\Bbb H}\cup {\Bbb Q}\cup \{\infty\}$  is 
the extended Lobachevsky plane [Gunning  1962]   \cite{Gun}. 
According to [Moeckel  1982]   \cite{Moe} {\it et al.}  such a choice of $S$
is the most natural for applications in number theory, e.g. the continued 
fractions.  For the sake of clarity,   we let $G\cong\Gamma_0(N)$, 
where   $N\ge 1$ is an integer  and 
$\Gamma_0(N)=\{(a,b,c,d)\in SL(2,{\Bbb Z}) ~|~c\equiv 0~mod ~N\}$ is 
the Hecke subgroup of $SL(2,{\Bbb Z})$. 
The modular surface $X_0(N)={\Bbb H}^*/\Gamma_0(N)$ is a 
complete hyperbolic surface.

Recall, that an {\it axis} is a closed geodesic {\bf g} $\in S$, which is  covered by a half-circle $\tilde{\hbox{\bf g}}$
$\in{\Bbb H}^*$  fixed by a hyperbolic transformation, see  [Casson \&  Bleiler   1988]   \cite{CaB} and [Gunning  1962]   \cite{Gun}. 
The lamination is called a {\it singleton},
if it consists of a unique simple closed geodesic {\bf g}; the singletons 
are dense in the space $\Lambda(S)$   [Canary, Epstein \&  Green   1987]   \cite{CEG}, Lemma 4.2.15.
The type (iii) lamination $\lambda\in\Lambda(X_0(N))$ will be called
a {\it Legendre lamination}, if there exists a regular continued 
fraction $[p_0,p_1, \dots]$, such that for $k\ge 1$ the axis of the hyperbolic
transformation 
\begin{equation}\label{eq1}
\left(\matrix{1 & p_0\cr 0 & 1}\right)
\left(\matrix{0 & 1 \cr 1 & p_1}\right)
\dots
\left(\matrix{0 & 1 \cr 1 & p_{2k}}\right)
\end{equation}
is a singleton {\bf g}$_k$ and $\lambda=\lim_{k\to\infty}$  {\bf g}$_k$.
We show in Lemma \ref{lm1} that the Legendre laminations exist and form an uncountable 
subset of $\Lambda(X_0(N))$.  The slope of the Legendre lamination is defined as 
\begin{equation}
\theta=p_0+{1\over {p_1+{1\over p_2+\dots}}}:=[p_0, p_1, p_2,\dots]
\end{equation}

\bigskip
To recover a lamination from $\theta$,  one needs to specify
the number and type of the boundary components of the lamination [Casson \&  Bleiler   1988]   \cite{CaB}.  
Let $g=g(N)$ be the genus of the modular surface and $\Delta=(k_1,\dots,k_m)$
a finite set of the positive integers and half-integers,  such that 
$\sum k_i=2g-2$. The $\Delta$ is called  a {\it singularity data} of the
lamination;  see Section 2.3 for the details.  Our main result is  the following existence and rigidity theorem.
\begin{thm}\label{thm1}
For each $N\ge 1$ there exists a continuum of the Legendre laminations
$\lambda\in\Lambda(X_0(N))$;  their  topological  types  are
bijective  with the pairs $(\Theta,\Delta)$, 
where $\Theta\equiv \theta~mod~GL(2,{\Bbb Z})$
is the equivalence class of irrationals modulo the action
of the matrix group $GL(2,{\Bbb Z})$ and $\Delta$ the singularity
data of $\lambda$.
\end{thm}
The article is organized as follows. In 
Section 2 the geodesic laminations are reviewed.
Theorem \ref{thm1} is proved in Section 3.

\section{Geodesic laminations}
This section is a brief review of the geodesic laminations following [Casson \&  Bleiler   1988]   \cite{CaB}
and contains no original results. We refer the reader to the above cited monograph 
for  a  detailed account.

\subsection{The Chabauty topology}
Let $S$ be a connected complete hyperbolic surface. By a {\it geodesic} {\bf g}
$\in S$ we understand the maximal arc consisting of the locally shortest sub-arcs.
The geodesic {\bf g} is called {\it simple} if it has no self-crossing
points.  (We leave aside the non-trivial question of  existence of simple geodesics;
such geodesics make an uncountable set on any hyperbolic surface, albeit of measure
zero [Artin  1924]   \cite{Art}.)  
A  {\it geodesic lamination} on $S$ is a closed subset $\lambda$ of $S$,
which is a disjoint union of the simple geodesics. The geodesics
are called the {\it leaves} of $\lambda$. The lamination
$\lambda$ is called {\it minimal} if no proper subset of $\lambda$
is a geodesic lamination. The following lemma gives an idea of the
minimal laminations (our main object of study). 
\begin{prp}\label{prp1}
The minimal lamination $\lambda$ in a complete hyperbolic
surface $S$ is either a singleton (simple closed geodesic) or an
uncountable nowhere dense (Cantor) subset of $S$. 
\end{prp}
{\it Proof.} See Lemma 4.2.2 of [Canary, Epstein \&  Green   1987]   \cite{CEG} and Lemma 3.3 of
[Casson \&  Bleiler   1988]   \cite{CaB}.
$\square$

\medskip
The set of all minimal laminations on the surface $S$ will be
denoted by $\Lambda(S)$.
The laminations $\lambda,\lambda'$ on $S$ are said to be {\it topologically
conjugate} if there exists a homeomorphism $\varphi:S\to S$, 
such that each leaf of $\lambda$ through the point $x\in S$ goes
to the leaf of $\lambda'$ through the point $\varphi(x)$. Clearly,
the set  $\Lambda(S)$  splits into the equivalence classes of 
topological conjugacy (or, topological {\it types}).

Since $S$ is a complete hyperbolic surface, its universal
cover is the unit disk $D$. Any $\lambda\in\Lambda(S)$
lifts to a lamination $\tilde\lambda$ on $D$, which is
invariant under the action of the covering transformations.
Every leaf of $\tilde\lambda$ is given by an unordered pair
of points at the boundary of $D$, and therefore the space
of geodesics is homeomorphic to the  M\"obius band,
$M$.

Let $C(M)$ be the set of all closed subsets of $M$. The
{\it Chabauty topology}
\footnote{It is a tradition to reserve the term Chabauty 
topology for a topology on the set of all closed subgroups
of a locally compact group. However, Canary, Epstein and Green [Canary, Epstein \&  Green   1987]   \cite{CEG}
differ from the convention but point out that on metrizable  space
this is nothing but the Hausdorff topology.}
 on $C(M)$ is given by the Hausdorff
distance $d(X,Y)\le\varepsilon$ iff  $X\subseteq N_{\varepsilon}(Y)$
and $Y\subseteq N_{\varepsilon}(X)$, 
where $N_{\varepsilon}(X)$ ($N_{\varepsilon}(Y)$) is a 
$\varepsilon$-neighbourhood of the closed set $X\in C(M)$
($Y\in C(M)$). The function $d$ turns the set $\tilde\Lambda$ of the
laminations in $D$ into a compact metrizable Hausdorff space.
The Chabauty topology on $\Lambda(S)$ can be defined as a
factor-topology of the topology on $\tilde\Lambda$ under the
covering map.  Everywhere in below, the standard topology on the set $\Lambda(S)$
will be the Chabauty topology. The following statement, mentioned
in the introduction, will be critical.
\begin{prp}\label{prp2}
The subspace $\Lambda^*(S)$ made of the singletons (i.e. the 
simple closed geodesics in $S$) is
dense in the space $\Lambda(S)$.
\end{prp}
{\it Proof.} See [Canary, Epstein \&  Green   1987]   \cite{CEG}, Lemma 4.2.15.
$\square$

\subsection{The principal regions and boundary leaves}
If $\lambda\in\Lambda(S)$, then a component of $S-\lambda$ is called
a {\it principal (complementary) region} for $\lambda$. 
(Note that $S-\lambda$ may have several connected components.)
The leaves of $\lambda$, which form the boundary of a principal
region,  are called the {\it boundary leaves.} 
If $\lambda$ is minimal (which we always assume to be the case), then each boundary 
leaf is a dense leaf of $\lambda$, isolated from one side.   
Note that by the proposition \ref{prp1} the area of a minimal lamination is zero, hence 
$Area~(S-\lambda)=Area~S$ and the principal region is a complete hyperbolic surface
of the area $-2\pi\chi(S)$, where $\chi(S)=2-2g$ is the Euler characteristic
of the surface $S$.  
(More precisely, if one  takes  the closure of the lift and quotient by the isometry group,  then  one gets
a hyperbolic surface with  geodesic boundary.) 
If $U$ is a component of the preimage of
the principal region in $D$, then $U$ is a union of 
the ideal polygons $U_i$ in $D$ (see Fig.1)

\begin{figure}
\begin{picture}(400,100)(0,0)


\put(100,65){\circle{40}}
\put(300,65){\circle{40}}

\qbezier(115,76)(100,65)(115,54)
\qbezier(85,76)(100,65)(85,54)
\qbezier(90,82)(100,65)(110,82)
\qbezier(90,48)(100,65)(110,48)

\qbezier(115,76)(113,69)(118,70)
\qbezier(118,70)(110,65)(115,54)

\qbezier(85,76)(87,69)(82,70)
\qbezier(82,70)(90,65)(85,54)

\qbezier(90,48)(97,56)(100,45)
\qbezier(100,45)(103,56)(110,48)

\qbezier(90,82)(97,74)(100,85)
\qbezier(100,85)(103,74)(110,82)


\qbezier(320,65)(300,65)(310,80)
\qbezier(310,80)(300,65)(290,80)
\qbezier(290,80)(300,65)(280,65)
\qbezier(280,65)(300,65)(290,50)
\qbezier(290,50)(300,65)(310,50)
\qbezier(310,50)(300,65)(320,65)

\put(30,20){{\sf (i) Maximal number of ideal}}
\put(85,7){{\sf polygons}}
\put(230,20){{\sf (ii) Minimal number of ideal}}
\put(280,7){{\sf polygons}}

\end{picture}
\caption{The region $U=\sqcup ~U_i$ for the surface of genus $2$}
\end{figure}

\medskip\noindent
The hyperbolic area of an ideal $n_i$-gon $U_i$ is equal to $(n_i-2)\pi$ 
 [Casson \&  Bleiler   1988]   \cite{CaB}. Since $\sum Area~U_i=(4g-4)\pi$,
the number of the ideal polygons in $U$ is  finite.

\subsection{The singularity data}
To capture combinatorial structure of the lamination $\lambda\in\Lambda(S)$,
we shall need the following collection of data.  
Denote by $U_i^{(n_i)}$ an $i$-th ideal $n_i$-gon in the principal region
of the lamination $\lambda$; here $n_i\ge 3$ is an integer. 
It is known, that  $Area~U_i^{(n_i)}=(n_i-2)\pi$. On the other hand,
the total area of all ideal polygons must be equal to the hyperbolic
area of the surface $S$, i.e.  $\sum_{i=1}^m~Area~U_i^{(n_i)}=(4g-4)\pi$. 
Thus, one arrives at the equation
\begin{equation}
\sum_{i=1}^m{n_i-2\over 2}=2g-2.
\end{equation}
For simplicity, we let $k_i={1\over 2}(n_i-2)$; since $n_i\ge 3$,
the numbers  $k_i$ take integer and half-integer positive values. In the above
notation, our formula becomes $\sum_{i=1}^mk_i=2g-2$.

Let $\{U_1^{(n_1)},\dots, U_{m}^{(n_m)}\}$ be a collection of the ideal
polygons in the principal region of the lamination $\lambda\in\Lambda(S)$; 
let $\{k_1,\dots,k_m\}$ be the corresponding collection of $k_i$. 
The unordered  tuple $\Delta=(k_1,\dots,k_m)$ will be called a 
{\it singularity data}  of the lamination $\lambda$. 
For example, the singularity data of laminations with the principal regions  shown in
 Fig.1 (i) and (ii)  are $\Delta_1=({1\over 2}, {1\over 2}, {1\over 2}, {1\over 2})$
and $\Delta_2=(2)$, respectively.

Conversely, given a set of positive integers and half-integers,
such that $\sum_{i=1}^mk_i=2g-2$, there exists a lamination $\lambda$
on surface $S$ of genus $g$, which realizes the singularity data $(k_1,\dots, k_m)$;
this fact follows from [Hubbard \& Masur 1979] \cite{HuMa1}.

Finally, the term ``singularity'' is justified by the fact, that minimal
geodesic laminations are bijective with the measured foliations on the same
surface [Thurston 1997]   \cite{T}; under the bijection each ideal polygon $U_i^{(n_i)}$ 
corresponds to a singular point of index $-k_i$ of the foliation.

\section{Proof of theorem 1}
The proof is arranged into a series of lemmas, starting with the following elementary
\begin{lem}\label{lm2}
The continued fraction of the Legendre lamination is unique. 
\end{lem}
{\it Proof.} 
Recall a classical bijection between the rationals $r_k$ and 
finite continued fractions $[p_0,\dots,p_{2k}]$,
given by the formula:
\begin{equation}\label{eq4}
r_k=
p_0+
{1\over\displaystyle p_1+
{\strut 1\over\displaystyle p_2+\dots+
{\strut 1\over\displaystyle p_{2k}}}},
\end{equation}
see [Perron 1954]   \cite{P}; the bijection extends to the singletons, since each
$[p_0,\dots,p_{2k}]$ defines a singleton {\bf g}$_k$. Moreover,  such a bijection
extends  to the half-circles $\tilde{\hbox{{\bf g}}}_k\in {\Bbb H}^*$ which lie in the
same orbit of the group $\Gamma_0(N)$. We refer the reader to [Artin 1924] \cite{Art}.

To the contrary, let $\tilde\lambda$ be the preimage of a Legendre lamination on ${\Bbb H}^*$,  
such that 
$\tilde\lambda=\lim_{k\to\infty} \tilde{\hbox{{\bf g}}}_k=\lim_{k\to\infty}  \tilde{\hbox{{\bf g}}}_k'$,
where  $\tilde{\hbox{{\bf g}}}_k\ne  \tilde{\hbox{{\bf g}}}_k'$.
By the Artin bijection, one gets
a pair of the regular continued fractions convergent to the same limit.
It is known to be false,  see e.g. [Perron 1954]   \cite{P},  Satz 2.6,  p. 33. 
$\square$

\begin{lem}\label{lm1}
The Legendre laminations form an uncountable  subset of $\Lambda(X_0(N))$. 
\end{lem}
{\it Proof.} 
(i) The set of the Legendre laminations is non-empty. 
Indeed, let $tr: \Gamma_0(N)\to {\Bbb Z}$ be the trace function.
We denote by ${\Bbb Z}_0$ the subset of $tr(\Gamma_0(N))$ consisting of traces of the hyperbolic
transformations,  whose axes are singletons.
It is known, that  ${\Bbb Z}_0$  is an infinite set;  moreover, each arithmetic progression
contains an infinite number of elements of ${\Bbb Z}_0$,  e.g. [Birman \& Series   1984]   \cite{BiSe1}.

Denote by $\gamma_k\in\Gamma_0(N)$ the product (\ref{eq1}).
The matrix multiplication (in the first few terms) gives us:
\begin{eqnarray}\label{eq3}
tr(\gamma_0) &=& 2\cr
tr(\gamma_1) &=& p_0+p_1\cr
tr(\gamma_2) &=& 2+p_0p_1+p_1p_2\cr
tr(\gamma_3) &=& p_0+p_1+p_2+p_3+p_0p_1p_2+p_1p_2p_3,
\end{eqnarray}
where $p_0\in {\Bbb N}\cup\{0\}$ and $p_i\in {\Bbb N}$. 
It is clear, that one can find $p_0,p_1,p_2$, so that 
$tr(\gamma_1)$ and $tr(\gamma_2)$ belong to ${\Bbb Z}_0$;
indeed, it suffices to take $p_0=z_1-1, p_1=1$ and $p_2=z_2-z_1-1$
for any two points $z_1,z_2\in {\Bbb Z}_0$, such that $2\le z_1\le z_2-2$.
In general, assume (by induction) that $p_0,p_1,\dots, p_n$ satisfy
the condition $tr(\gamma_i)\in {\Bbb Z}_0$ for all $0\le i\le n$.
We want to choose $p_{n+1}$, such that the matrix $\gamma_{n+1}$ defined 
by (\ref{eq1}) is hyperbolic.  Consider an arithmetic progression
$ap_{n+1}+b=tr(\gamma_{n+1})$,
where $a$ and $b$ are integers depending only on $p_0,\dots,p_n$. The arithmetic
progression $\{a+b, 2a+b,\dots\}$ contains infinitely many
elements of ${\Bbb Z}_0$;  let  $z_{n+1}$ be one of them
with the index $p_{n+1}$. Then, $tr(\gamma_{n+1})\in {\Bbb Z}_0$
and the induction is completed.

In this way, one obtains an infinite sequence $[p_0,p_1,\dots]$,
such that the axes of hyperbolic transformations $\gamma_i$
are singletons; their hyperbolic length grows and, by compactness
of $\Lambda(X_0(N))$, converges to a type (iii) lamination $\lambda$.
By construction, $\lambda$ is the Legendre lamination.

\medskip
(ii) The Legendre laminations are uncountable; indeed, 
by a modification of the argument of (i), one can choose, at each step
of the induction, infinitely many  $p_{n+1}$,
which satisfy the equation  $ap_{n+1}+b=tr(\gamma_{n+1})$.
Thus, one gets an infinite sequence  $[p_0,p_1,\dots]$,
where each $p_i$ runs a countable infinite set; the set of 
all such sequences is also infinite, but uncountable.

In view of lemma \ref{lm2},   the corresponding Legendre laminations are 
also  uncountable;    lemma \ref{lm1} follows. 
 $\square$

\begin{lem}\label{lm3}
Two Legendre laminations are topologically conjugate,
if and only if their singularity data coincide and their continued fractions coincide,  except a finite
number of terms. 
\end{lem}
{\it Proof.}
Let $\varphi: X_0(N)\to X_0(N)$ be an automorphism, which conjugates
laminations $\lambda$ and $\lambda'=\varphi(\lambda)$.
Notice that the automorphism $\varphi$ preserves the singularity data
of $\lambda$ and $\lambda'$,  see Section 2.3.    
The action of $\varphi$ extends to the singletons; by the Artin bijection,
it extends to the positive rationals $r_k$,  cf.  proof of lemma \ref{lm2}. 
But each automorphism of the rational numbers is given by the formula:
\begin{equation}\label{eq5}
\varphi(r_k)={ar_k+b\over cr_k+d},\quad
\left(\matrix{a & b\cr c & d}\right)\in GL(2,{\Bbb Z}).
\end{equation}
In the last formula, we can pass to the limit
$\theta=\lim_{k\to\infty}r_k$; thus, whenever
the $\theta$ is given by the continued fraction $[p_0,p_1,\dots]$,
then $\theta'=\varphi(\theta)$ is given by continued fraction 
 $[q_1,\dots,q_k; p_0,p_1,\dots]$, where 
\begin{equation}\label{eq6}
\left(\matrix{0 & 1\cr 1 & q_1}\right)\dots
\left(\matrix{0 & 1\cr 1 & q_k}\right)=
\left(\matrix{a & b\cr c & d}\right).
\end{equation}
Notice, that the Legendre laminations $\lambda$ and $\lambda'$ 
are given by the fractions $\theta$ and $\theta'$, respectively.
The converse statement can be proved similarly. 
Lemma \ref{lm3} follows.
$\square$

\begin{lem}\label{lm4}
For each $\theta\in {\Bbb R}-{\Bbb Q}$ and an abstract data $\Delta$
there exists a Legendre lamination of slope $\theta$,  which realizes 
the   singularity data $\Delta$.  
\end{lem}
{\it Proof.}
Let $D$ be the unit disk and $\Delta=(k_1,\dots,k_m)$ a singularity data compatible with the modular surface $X_0(N)$;
denote by $U_1,\dots,U_m$
the ideal polygons corresponding to $\Delta$. 
(Notice that the polygons $U_i$ are not uniquely defined;  yet their isometry class in the
group $\Gamma_0(N)$ depends solely on $\Delta$.)
If $\theta=[p_0,p_1,\dots]$,  we shall write $C_i$ to denote an axis of transformation 
$\gamma_i$ computed by formula (\ref{eq1});   we denote the corresponding singleton
by $\lambda_i^*$.   For each $\gamma_i$ consider a region ${\cal O}_{\gamma_i}$
of $D$ defined by the formula:  
\begin{equation}\label{eq7}
{\cal O}_{\gamma_i}=\bigcup_{g\in\Gamma_0(N)}
~\bigcup_{|tr(h)|\le tr(\gamma_i)}
ghg^{-1} 
\left(\sum_{i=1}^mU_i\right);
\end{equation}
notice that region  ${\cal O}_{\gamma_i}$ is invariant of the conjugacy
class of transformation $\gamma_i$.  We denote the maximal subset of  $D-{\cal O}_{\gamma_i}$
containing $C_i$ (and all  isometric images of  $C_i$) by $\Omega_i$.

Let us show  that $D\supset\Omega_0\supset\Omega_1\supset\dots$ are strict inclusions.
Indeed,  in view of formulas (\ref{eq3}),   we have $tr~(\gamma_{i+1})>tr~(\gamma_i)$
for all $i\ge0$;  thus for the cosets $H_i:=\{ghg^{-1} ~|~ g\in\Gamma_0(N), ~|tr (h)|\le tr(\gamma_i)\}$, 
one gets an inclusion $H_i\subset H_{i+1}$.   In view of  formula (\ref{eq7}),  one obtains an inclusion
${\cal O}_{\gamma_i}\subset {\cal O}_{\gamma_{i+1}}$;   thus,  for the complement sets $\Omega_i$, we
have  a strict  inclusion $\Omega_i\supset \Omega_{i+1}$,  which holds for all $i\ge0$.

Denote by $\Omega=\cap_{i=0}^{\infty}\Omega_i$. It is easy to see,  that
 $\Omega$ is a non empty closed set as limit of a decreasing sequence of non empty closed sets.

By our construction, each $\Omega_i$  contains the arc $C_i$;   therefore, the set 
$C=\lim_{i\to\infty} C_i$ is contained in the set $\Omega$.  But $\Omega$ is a closed
set  (of measure zero)  of the unit disk $D$;   thus $C$ covers a type (iii) lamination 
$\lambda=\lim_{i\to\infty}$ {\bf g}$_i$   on the surface $X_0(N)$.  The lamination $\lambda$ is
a Legendre lamination of slope $\theta$,  which realizes the singularity data $\Delta$.
Lemma \ref{lm4} is proved.  
$\square$

\bigskip
Now we use the above four lemmas to prove theorem \ref{thm1}. 
Namely,  let $\Lambda_0\subset\Lambda(X_0(N))$ be the space of all Legendre
lamination on $X_0(N)$; in view of lemmas \ref{lm2} and \ref{lm3}, each conjugacy class of 
$\lambda\in\Lambda_0$ defines an invariant $(\Theta,\Delta)$,  where 
$\Theta=\{\theta' : \theta'={a\theta+b\over c\theta +d}; ~ a,b,c,d\in {\Bbb Z}, ~ad-bc=\pm 1\}$.
Conversely,   lemma \ref{lm4}  says that any  abstractly given  invariant  $(\Theta,\Delta)$
admits a  realization by a  lamination $\lambda\in\Lambda_0$. 
This argument finishes the proof of theorem \ref{thm1}.    
$\square$

\bigskip\noindent
{\sf Acknowledgements.} I am grateful to Samuil ~Kh. ~Aranson  and Evgeny ~V. ~Zhuzhoma for an introduction
to  Weil's problem.    I thank the referee for a proofreading of the manuscript.



\vskip1cm

\textsc{Department of Mathematics and Computer Science, St.~John's University, 8000 Utopia Parkway,  
New York,  NY 11439, United States;} ~\textsc{E-mail:} {\sf igor.v.nikolaev@gmail.com}

\end{document}